\documentclass[12pt]{article}
\usepackage{graphicx}
\usepackage{amsmath,amsthm,amssymb,enumerate}
\usepackage{euscript,mathrsfs}
\usepackage{color}
\usepackage{dsfont}
\usepackage[left=2cm,right=2cm,top=3.5cm,bottom=3.5cm]{geometry}
\usepackage{color}
\usepackage[framemethod=tikz]{mdframed}
\allowdisplaybreaks

\usepackage{soul}

\catcode`\@=11 \@addtoreset{equation}{section}

\catcode`\@=12

\allowdisplaybreaks

\newtheorem{Theorem}{Theorem}[section]
\newtheorem{Proposition}[Theorem]{Proposition}
\newtheorem{Lemma}[Theorem]{Lemma}
\newtheorem{Corollary}[Theorem]{Corollary}

\theoremstyle{definition}
\newtheorem{Definition}[Theorem]{Definition}

\newtheorem{Remark}[Theorem]{Remark}

\newcommand{\bTheorem}[1]{
\begin{Theorem} \label{T#1} }
\newcommand{\eT}{\end{Theorem}}

\newcommand{\bProposition}[1]{
\begin{Proposition} \label{P#1}}
\newcommand{\eP}{\end{Proposition}}

\newcommand{\bLemma}[1]{
\begin{Lemma} \label{L#1} }
\newcommand{\eL}{\end{Lemma}}

\newcommand{\bCorollary}[1]{
\begin{Corollary} \label{C#1} }
\newcommand{\eC}{\end{Corollary}}

\newcommand{\bRemark}[1]{
\begin{Remark} \label{R#1} }
\newcommand{\eR}{\end{Remark}}

\newcommand{\bDefinition}[1]{
\begin{Definition} \label{D#1} }
\newcommand{\eD}{\end{Definition}}

\newcommand{\Del}{\Delta_x}

\newcommand{\Ds}{\mathbb{D}_x}

\newcommand{\bfomega}{\boldsymbol{\omega}}

\newcommand{\bfphi}{\boldsymbol{\varphi}}

\newcommand{\bFormula}[1]{
\begin{equation} \label{#1}}
\newcommand{\eF}{\end{equation}}

\newcommand{\Ov}[1]{\overline{#1}}

\newcommand{\aleq}{\stackrel{<}{\sim}}

\newcommand{\vr}{\varrho}
\newcommand{\vre}{\vr_\ep}

\newcommand{\vte}{\vt_\ep}
\newcommand{\vue}{\vu_\ep}

\newcommand{\tvt}{\tilde \vt}

\newcommand{\vt}{\vartheta}
\newcommand{\vu}{\vc{u}}

\newcommand{\vc}[1]{{\bf #1}}

\newcommand{\Div}{{\rm div}_x}
\newcommand{\Grad}{\nabla_x}

\newcommand{\dx}{\,{\rm d} {x}}

\newcommand{\dt}{\,{\rm d} t }

\newcommand{\intO}[1]{\int_{\Omega} #1 \ \dx}

\newcommand{\intST}[1]{\int_{S_T} #1 \ \dt}

\newcommand{\D}{{\rm d}}

\newcommand{\ep}{\varepsilon}

\newcommand{\vtB}{\vt_B}

\newcommand{\br}{ \nonumber \\ }

\def\softd{{\leavevmode\setbox1=\hbox{d}%
          \hbox to 1.05\wd1{d\kern-0.4ex{\char039}\hss}}}
\definecolor{Cgrey}{rgb}{0.85,0.85,0.85}
\definecolor{Cblue}{rgb}{0.50,0.85,0.85}
\definecolor{Cred}{rgb}{1,0,0}
\definecolor{fancy}{rgb}{0.10,0.85,0.10}

\newcommand\Cbox[2]{%
    \newbox\contentbox%
    \newbox\bkgdbox%
    \setbox\contentbox\hbox to \hsize{%
        \vtop{
            \kern\columnsep
            \hbox to \hsize{%
                \kern\columnsep%
                \advance\hsize by -2\columnsep%
                \setlength{\textwidth}{\hsize}%
                \vbox{
                    \parskip=\baselineskip
                    \parindent=0bp
                    #2
                }%
                \kern\columnsep%
            }%
            \kern\columnsep%
        }%
    }%
    \setbox\bkgdbox\vbox{
        \color{#1}
        \hrule width  \wd\contentbox %
               height \ht\contentbox %
               depth  \dp\contentbox
        \color{black}
    }%
    \wd\bkgdbox=0bp%
    \vbox{\hbox to \hsize{\box\bkgdbox\box\contentbox}}%
    \vskip\baselineskip%
}

\mdfdefinestyle{MyFrame}{%
	linecolor=black,
	outerlinewidth=1pt,
	roundcorner=5pt,
	innertopmargin=\baselineskip,
	innerbottommargin=\baselineskip,
	innerrightmargin=10pt,
	innerleftmargin=10pt,
	backgroundcolor=white!20!white}

\date{}

\begin{document}


\title{Time periodic motion of temperature driven compressible fluids}

\author{Eduard Feireisl
\thanks{The work of E.F. was partially supported by the
Czech Sciences Foundation (GA\v CR), Grant Agreement
21--02411S. The Institute of Mathematics of the Academy of Sciences of
the Czech Republic is supported by RVO:67985840. This work is partially supported by the Simons Foundation Award 
No 663281 granted to the Institute of Mathematics of the Polish Academy of Sciences for the years 2021-2023.} \and Piotr Gwiazda \and Agnieszka \'Swierczewska-Gwiazda
\thanks{The work of A. \'S-G. was partially supported by  National Science Centre
(Poland),  agreement no 2017/27/B/ST1/01569.}}

\date{}

\maketitle

\medskip

\centerline{Institute of Mathematics of the Academy of Sciences of the Czech Republic}

\centerline{\v Zitn\' a 25, CZ-115 67 Praha 1, Czech Republic}

\medskip

\centerline{Institute of Applied Mathematics and Mechanics, University of Warsaw}

\centerline{Banacha 2, 02-097 Warsaw, Poland}

\begin{abstract}
	
We consider the Navier--Stokes--Fourier system describing the motion of a compressible viscous fluid in a container with impermeable boundary subject to time periodic 
heating and under the action of a time periodic potential force. We show the existence of a time periodic weak solution for arbitrarily large physically admissible data.

\end{abstract}




\tableofcontents

\section{Introduction}
\label{i}

There are numerous examples of turbulent fluid motion excited by changes of the boundary temperature, among which the well studied problem of Rayleigh--B\' enard convection, see see e.g. Davidson \cite{DAVI}.
Motivated by similar problems in astrophysics of gaseous stars, we consider a general compressible viscous possibly rotating fluid, occupying a bounded domain $\Omega \subset R^d$, $d=2,3$, driven 
by periodic changes of boundary temperature. The relevant system of field equations for the standard variables: the mass density $\vr = \vr(t,x)$, the velocity $\vu = \vu(t,x)$, and the 
(absolute) temperature $\vt = \vt(t,x)$ reads:

\begin{mdframed}[style=MyFrame]
\begin{align}
	\partial_t \vr + \Div (\vr \vu) &= 0, \label{i1} \\
	\partial_t (\vr \vu) + \Div (\vr \vu \otimes \vu) + \vr (\boldsymbol{\omega} \times \vu) + \Grad p (\vr, \vt) &= \Div \mathbb{S} + \vr \Grad G, \label{i2}\\
	\partial_t (\vr e(\vr, \vt)) + \Div (\vr e(\vr, \vt) \vu) + \Grad \vc{q} &= \mathbb{S} : \Ds \vu - p(\vr, \vt) \Div \vu,
	\label{i3}
	\end{align}
\end{mdframed}

\noindent where $\mathbb{S}$ is the viscous stress given by Newton's rheological law
\begin{equation} \label{i4}
	\mathbb{S} (\vt, \Ds \vu) = \mu(\vt) \left( \Grad \vu + \Grad^t \vu - \frac{2}{3} \Div \vu \mathbb{I} \right) + 
	\eta(\vt) \Div \vu \mathbb{I},
\end{equation}
and $\vc{q}$ is the heat flux given by  
Fourier's law 
\begin{equation} \label{i5}
	\vc{q}(\vt, \Grad \vt) = - \kappa(\vt) \Grad \vt.
\end{equation} 
The momentum equation is augmented by the Coriolis force with the rotation  constant vector $\bfomega$, the associated centrifugal force as well as the gravitation and other possible inertial time-periodic forces are regrouped 
in the potential $G$. The fluid occupies a bounded smooth domain 
$\Omega \subset R^d$, $d=2,3$ endowed with the Dirichlet boundary conditions
\begin{mdframed}[style=MyFrame]
	
	\begin{align}
	\vu|_{\partial \Omega} = 0, \br 
	\vt|_{\partial \Omega} = \vtB \label{i6}
\end{align}	
\end{mdframed} 	
The functions $\vtB =\vtB(t,x)$ and $G = G(t,x)$ are smooth and $T-$periodic in the time variable,
\begin{mdframed}[style=MyFrame]
	
	\begin{align}
		\vtB(t + T,x) = \vtB(t,x), \br 
		G(t + T,x) = G(t,x). \label{i7}
	\end{align}	
\end{mdframed}
Hereafter, the problem \eqref{i1}--\eqref{i6} is referred to as \emph{Navier--Stokes--Fourier system}.

Our goal is to show the existence of a time--periodic solution to problem \eqref{i1}--\eqref{i7}. There is a substantial number of references, where such a result is proven under some smallness and smoothness 
assumption on the data. Valli and Zajaczkowski \cite{Valli1}, \cite{VAZA} observe that the distance of two smooth global in time solutions decays in time for the system close to a stable equilibrium, and, as a by product, they deduce the existence of a time periodic solution. Similar ideas have been followed by many authors, see B\v rezina and Kagei \cite{BreKag2}, \cite{BreKag1}, Jin and Yang \cite{JinYang} , Kagei and Oomachi \cite{KagOom}, 
Kagei and Tsuda \cite{KagTsu}, Tsuda \cite{Tsuda} to name only a few. Turbulent fluid flows given by large forces out of equilibrium are mostly considered in the framework of \emph{weak} solutions. 
Based on the mathematical theory of compressible fluids developed by Lions \cite{LI}, \cite{LI4}, the existence of large time periodic solutions for the simplified isentropic system was proved 
in \cite{FMPS} for the isentropic pressure--density equation of state $p(\vr) = a \vr^\gamma$, $\gamma \geq \frac{9}{5}$. The later development of the theory in \cite{FNP} enabled to extend the result 
to the case $\gamma > \frac{5}{3}$, see Cai and Tan \cite{CaiTan}.

The situation is more delicate for the complete fluid systems including thermal effects. As a direct consequence of the Second law of thermodynamics, the existence of (forced) time periodic solutions 
is ruled out for problems with purely conservative boundary conditions, see \cite{FP20}. In \cite{FeMuNoPo}, the heat flux was controlled by means of a Robin type boundary condition 
\begin{equation} \label{i8}
	\vc{q} \cdot \vc{n} = d (\vt - \Theta_0) \ \mbox{on}\ \partial \Omega,
\end{equation}
with a given ``mean'' temperature $\Theta_0$. Accordingly, the internal energy is transferred out of the fluid domain in the high temperature regime and the time periodic motion is possible, see 
\cite[Theorem 1]{FeMuNoPo}. Our goal is to show a similar result for the Dirichlet boundary conditions \eqref{i6}. Note that the problem is much more delicate than in \cite{FeMuNoPo} as the heat flux through the boundary is {\it a priori} not controlled.

Our approach is based on several rather new ideas that appeared only recently in the mathematical theory of open fluid systems.

\begin{itemize}
	
	\item The concept of \emph{weak solution} for the Navier--Stokes--Fourier system based on a combination of the entropy inequality and the ballistic energy balance developed in \cite{ChauFei}.
	
	\item \emph{Uniform bounds} and large time asymptotics of the weak solutions in the spirit of \cite{FeiSwGw}.
	
	\item An \emph{approximation scheme} based on a penalization of the Dirichlet boundary conditions via \eqref{i8}.
	
	\end{itemize}

The constitutive restrictions imposed on the equations of state as well as the transport coefficients are the same as in the existence theory \cite{ChauFei}. In particular, the general equation of state 
of real monoatomic gases proposed in \cite[Chapters 1,2]{FeNo6A} is included. From this point of view, the result is apparently better than in the isentropic case studied in \cite{FMPS}, and later revisited by 
Cai and Tan \cite{CaiTan}, where the condition $\gamma > \frac{5}{3}$ is needed. The price to pay is the potential form of the driving force $\vc{f} = \Grad G$ that, however, includes the physically relevant 
centrifugal as well as gravitational forces.

The paper is organized as follows. In Section \ref{M}, we introduce the basic hypotheses concerning the constitutive relations and state the main result. In Section \ref{A}, we introduce an approximation scheme 
inspired by \cite{FeMuNoPo}. Section \ref{U} is the heart of the paper. Here we establish the necessary uniform bounds to perform the limit in the sequence of approximate solutions. Finally, in Section \ref{C}, we obtain the desired solution as a limit of the approximate sequence.

\section{Main result}
\label{M}

Before stating the main result, we recall the form of the constitutive equations proposed in \cite[Chapters 1,2]{FeNo6A}. To comply with the Second law of thermodynamics, we postulate the existence of \emph{entropy} $s$, related 
to the internal energy $e$ and the pressure $p$ through Gibbs' equation
\begin{equation} \label{M1}
	\vt D s(\vr, \vt) = D e(\vr, \vt) + p(\vr, \vt) D \left( \frac{1}{\vr} \right).
	\end{equation} 
 
\subsection{Constitutive theory}
Similarly to \cite[Chapters 1,2]{FeNo6A} we consider the pressure equation of state in the form 
\[
p(\vr, \vt) = p_{\rm m} (\vr, \vt) + p_{\rm rad}(\vt), 
\]
where $p_{\rm m}$ is the pressure of a general \emph{monoatomic} gas related to the internal energy through
\begin{equation} \label{con1}
	p_{\rm m} (\vr, \vt) = \frac{2}{3} \vr e_{\rm m}(\vr, \vt),
\end{equation}
augmented by the radiation pressure 
\[
p_{\rm rad}(\vt) = \frac{a}{3} \vt^4,\ a > 0.
\]
Similarly, the internal energy reads 
\[
e(\vr, \vt) = e_{\rm m}(\vr, \vt) + e_{\rm rad}(\vr, \vt),\ e_{\rm rad}(\vr, \vt) = \frac{a}{\vr} \vt^4.
\]

Now, Gibbs' equation \eqref{M1} gives rise to a specific form of $p_m$, 
\[
p_m (\vr, \vt) = \vt^{\frac{5}{2}} P \left( \frac{\vr}{\vt^{\frac{3}{2}}  } \right)
\]
for a certain $P \in C^1[0,\infty)$.
	Consequently, 
	\begin{equation} \label{M2}
		p(\vr, \vt) = \vt^{\frac{5}{2}} P \left( \frac{\vr}{\vt^{\frac{3}{2}}  } \right) + \frac{a}{3} \vt^4,\ 
		e(\vr, \vt) = \frac{3}{2} \frac{\vt^{\frac{5}{2}} }{\vr} P \left( \frac{\vr}{\vt^{\frac{3}{2}}  } \right) + \frac{a}{\vr} \vt^4, \ a > 0.
	\end{equation}
In addition, we suppose	
\begin{equation} \label{w10}
		P(0) = 0,\ P'(Z) > 0 \ \mbox{for}\ Z \geq 0,\ 0 < \frac{ \frac{5}{3} P(Z) - P'(Z) Z }{Z} \leq c \ \mbox{for}\ Z > 0.
	\end{equation} 
that may be seen as a direct consequence of \emph{hypothesis of thermodynamic stability}, see \cite[Chapter 1]{FeNo6A}, and Bechtel, Rooney, and Forrest \cite{BEROFO}.
It follows that the function $Z \mapsto P(Z)/ Z^{\frac{5}{3}}$ is decreasing, and we suppose 
	\begin{equation} \label{w11}
		\lim_{Z \to \infty} \frac{ P(Z) }{Z^{\frac{5}{3}}} = p_\infty > 0.
	\end{equation}
	
	 The associated entropy
	takes the form 
	\begin{equation} \label{w12}
		s(\vr, \vt) = \mathcal{S} \left( \frac{\vr}{\vt^{\frac{3}{2}} } \right) + \frac{4a}{3} \frac{\vt^3}{\vr}, 
	\end{equation}
	where 
	\begin{equation} \label{w13}
		\mathcal{S}'(Z) = -\frac{3}{2} \frac{ \frac{5}{3} P(Z) - P'(Z) Z }{Z^2} < 0.
	\end{equation}
Finally, we impose the Third law of thermodynamics, see e.g. Belgiorno \cite{BEL1}, \cite{BEL2}, requiring the total entropy to vanish 
	as soon as the absolute temperature approaches zero, 
	\begin{equation} \label{w14}
		\lim_{Z \to \infty} \mathcal{S}(Z) = 0.
	\end{equation}
It is easy to check that \eqref{w10} -- \eqref{w14} imply
\begin{equation} \label{w15}
	0 \leq	\vr \mathcal{S} \left( \frac{\vr}{\vt^{\frac{3}{2}}} \right) \leq c \left(1 + \vr \log^+(\vr) + \vr \log^+(\vt) \right).
\end{equation}

As for the transport coefficients, we suppose that they are continuously differentiable functions of the absolute temperature satisfying
\begin{align} 
	0 < \underline{\mu}(1 + \vt) &\leq \mu(\vt),\ |\mu'(\vt)| \leq \Ov{\mu}, \br 
	0 \leq \eta (\vt) &\leq \Ov{\eta}(1 + \vt), \br
	0 < \underline{\kappa} (1 + \vt^\beta) &\leq \kappa (\vt) \leq \Ov{\kappa}(1 + \vt^\beta), \label{w16}
\end{align}
where, in accordance with the existence theory developed in \cite{ChauFei}, we require 
\begin{equation} \label{w18}
	\beta > 6.
\end{equation}

\subsection{Weak solutions}

It is convenient to identify the time periodic functions (distributions) 
with objects defined on a periodic ``flat sphere''
\[
S_T = [0,T]|_{\{ 0, T \} }
\]
We are ready to introduce the concept of \emph{time periodic solution} to the Navier--Stokes--Fourier system 
\eqref{i1}--\eqref{i7}.

\begin{mdframed}[style=MyFrame]
	
		\begin{Definition} [{\bf weak solution}] \label{Dw1}
	
	We say that a trio $(\vr, \vt, \vu)$ is a weak time--periodic solution to the problem \eqref{i1}--\eqref{i7}
	if the following holds:

	\begin{itemize}
		
		\item {\bf Regularity class}: 
		\begin{align}
			\vr &\in C_{\rm weak}(S_T; L^\gamma(\Omega)) \ \mbox{for}\ \gamma = \frac{5}{3}, \br
			\vu &\in L^2(S_T; W^{1,2}_0 (\Omega; R^d)),\ 
			\vr \vu \in C_{\rm weak}(S_T, L^{\frac{2 \gamma}{\gamma + 1}}(\Omega; R^d)), \br 
			\vt^{\beta/2} ,\ \log(\vt) &\in L^2(S_T; W^{1,2}(\Omega)) ,\br
			(\vt - \vtB) &\in L^2(S_T; W^{1,2}_0 (\Omega)).
			\label{w6}
			\end{align}
		 
		\item {\bf Equation of continuity:} 
	
		\begin{align} 
			\intST{ \intO{ \left[ \vr \partial_t \varphi + \vr \vu \cdot \Grad \varphi \right] } } &= 0, 
			\label{w3} \\
			\intST{ \intO{ \left[ b(\vr) \partial_t \varphi + b(\vr) \vu \cdot \Grad \varphi + \Big( 
				b(\vr) - b'(\vr) \vr \Big) \Div \vu \varphi \right] } } &=0
			\label{w4}
			\end{align}
for any $\varphi \in C^1(S_T \times \Ov{\Omega} )$, and any $b \in C^1(R)$, 	$b' \in C_c(R)$.
\item {\bf Momentum equation:} 
\begin{align}
&\intST{ \intO{ \Big[ \vr \vu \cdot \partial_t \bfphi + \vr \vu \otimes \vu : \Grad \bfphi - 
		\vr (\bfomega \times \vu) \cdot \bfphi + 
	p \Div \bfphi \Big] } } \br &= \intST{ \intO{ \Big[ \mathbb{S} : \Grad \bfphi - \vr \Grad G \cdot \bfphi \Big] } }, 
\label{w5}
\end{align}	
for any $\bfphi \in C^1_c(S_T \times \Omega; R^d)$.

\item {\bf Entropy inequality:}
\begin{align}
	- &\intST{ \intO{ \left[ \vr s \partial_t \varphi + \vr s \vu \cdot \Grad \varphi + \frac{\vc{q}}{\vt} \cdot 
		\Grad \varphi \right] } } \br &\geq \intST{ \intO{ \frac{\varphi}{\vt} \left[ \mathbb{S} : \Ds \vu - 
		\frac{\vc{q} \cdot \Grad \vt }{\vt} \right] } }
	\label{w7} 
	\end{align}
for any $\varphi \in C^1_c(S_T \times \Omega)$, $\varphi \geq 0$;

\item {\bf Ballistic energy balance:}
\begin{align}  
	- &\intST{ \partial_t \psi	\intO{ \left[ \frac{1}{2} \vr |\vu|^2 + \vr e - \tvt \vr s \right] } } + \intST{ \psi
	\intO{ \frac{\tvt}{\vt}	 \left[ \mathbb{S}: \Ds \vu - \frac{\vc{q} \cdot \Grad \vt }{\vt} \right] } }  \br
	&\leq 
	\intST{ \psi \intO{ \left[ \vr \vu \cdot \Grad G - \vr s \vu \cdot \Grad \tvt - \frac{\vc{q}}{\vt} \cdot \Grad \tvt \right] } }
	\label{w8}
\end{align}
for any $\psi \in C^1(S_T)$, $\psi \geq 0$, and any $\tvt \in C^1(S_T \times \Ov{\Omega})$, 
\begin{equation} \label{w9}
\tvt > 0,\ \tvt|_{\partial \Omega} = \vtB.
\end{equation}
		\end{itemize}
	
	\end{Definition}
	
	\end{mdframed}

The weak time--periodic solutions are therefore the weak solutions in the sense of \cite{ChauFei} that are 
$T-$periodic in the time variable. The instantaneous values the 
\emph{conservative variables} $\vr(\tau, \cdot)$, $(\vr \vu)(\tau, \cdot)$ are well defined as well as the 
right and left hand limits of the total entropy $S = \vr s(\vr, \vt)$, 
\begin{align}
\left< S(\tau-, \cdot); \phi \right> &\equiv \lim_{\delta \to 0+} \frac{1}{\delta} \int_{\tau - \delta}^\tau \intO{ \vr s (t,\cdot) \phi } \dt, \ 
\left< S(\tau+, \cdot); \phi \right> \equiv \lim_{\delta \to 0+} \frac{1}{\delta} \int_{\tau}^{\tau + \delta} \intO{ \vr s(t,\cdot) \phi } \dt.	
\nonumber
	\end{align}

\subsection{Main result}

Having collected the necessary preliminary material we are ready to state our main result. 

\begin{mdframed}[style=MyFrame]
	
	\begin{Theorem}[{\bf existence of time periodic solutions}] \label{MT1}
		Let $\Omega \subset R^d$, $d=2,3$ be a bounded domain of class $C^{2 + \nu}$. Suppose that the pressure 
		$p$, the internal energy $e$, the entropy $s$, as well as the transport coefficients $\mu$, $\eta$, and $\kappa$ satisfy the hypotheses \eqref{con1}--\eqref{w18}. Finally, let the data
		$G \in W^{1,\infty}(S_T \times \Omega)$, $\vtB \in C^3(S_T \times R^d)$ 
		be time periodic as stated in \eqref{i7}, and 
		\[
		\inf_{S_T \times \Omega} \vtB = \underline{\vt} > 0.
		\]
		
		Then for any $M_0$ there exists at least one time periodic solution $(\vr, \vt, \vu)$ of the problem \eqref{i1}--\eqref{i7} in the sense specified in Definition \ref{Dw1} satisfying 
		\[
		\intO{ \vr(t,\cdot) } = M_0 \ \mbox{for any}\ t \in S_T. 
		\]
		\end{Theorem}
	
	\end{mdframed}

\begin{Remark} \label{MR1} In the hypotheses of Theorem \ref{MT1}, we assume that $\vtB|_{\partial \Omega}$ is a restriction of a (smooth) function defined on the whole space $R^d$.
	\end{Remark}

The rest of the paper is devoted to the proof of Theorem \ref{MT1}.

\section{Approximate problem}

\label{A}

The most efficient way of constructing suitable approximate solutions seems adapting the result of 
\cite{FeMuNoPo} to the present setting. Specifically, the approximation scheme is based on penalization 
the Dirichlet boundary condition for the temperature via the Robin boundary conditions 
\begin{equation} \label{A1}
	\vc{q} \cdot \vc{n} = \frac{1}{\ep} |\vt - \vtB|^{k}(\vt - \vtB),\ k \geq 0, \ \mbox{on}\ \partial \Omega,
	\end{equation}
where $\ep > 0$ is a small parameter.

The approximate solutions $(\vre, \vte, \vue)$ are defined similarly to Definition \ref{Dw1}:

			\begin{itemize}
			
			\item {\bf Regularity class:} 
			\begin{align}
				\vre &\in C_{\rm weak}(S_T; L^\gamma(\Omega)) \ \mbox{for}\ \gamma = \frac{5}{3}, \br
				\vue &\in L^2(S_T; W^{1,2}_0 (\Omega; R^d)),\ 
				\vre \vue \in C_{\rm weak}(S_T, L^{\frac{2 \gamma}{\gamma + 1}}(\Omega; R^d)) \br 
				\vte^{\beta/2} ,\ \log(\vte) &\in L^2(S_T; W^{1,2}(\Omega)).
				\label{a1}
			\end{align}
			
			\item {\bf Equation of continuity:} 
			
			\begin{align} 
				\intST{ \intO{ \left[ \vre \partial_t \varphi + \vre \vue \cdot \Grad \varphi \right] } } &= 0, 
				\label{w3A} \\
				\intST{ \intO{ \left[ b(\vre) \partial_t \varphi + b(\vre) \vue \cdot \Grad \varphi + \Big( 
						b(\vre) - b'(\vre) \vre \Big) \Div \vue \varphi \right] } } &=0
				\label{a2}
			\end{align}
			for any $\varphi \in C^1(S_T \times \Ov{\Omega} )$, and any $b \in C^1(R)$, 	$b' \in C_c(R)$.
			\item {\bf Momentum equation:} 
			\begin{align}
				&\intST{ \intO{ \Big[ \vre \vue \cdot \partial_t \bfphi + \vre \vue \otimes \vue : \Grad \bfphi - 						\vre (\bfomega \times \vue) \cdot \bfphi + 
						p \Div \bfphi \Big] } } \br &= \intST{ \intO{ \Big[ \mathbb{S} : \Grad \bfphi - \vre \Grad G \cdot \bfphi \Big] } }, 
				\label{a3}
			\end{align}	
			for any $\bfphi \in C^1_c(S_T \times \Omega; R^d)$.
			
			\item {\bf Entropy inequality:}
			\begin{align}
				- &\intST{ \intO{ \left[ \vre s \partial_t \varphi + \vre s \vue \cdot \Grad \varphi + \frac{\vc{q}}{\vte} \cdot 
						\Grad \varphi \right] } } \br &\geq \intST{ \intO{ \frac{\varphi}{\vte} \left[ \mathbb{S} : \Ds \vue - 
						\frac{\vc{q} \cdot \Grad \vte }{\vte} \right] } } + \frac{1}{\ep} 
					\intST{ \int_{\partial \Omega} \varphi \frac{|\vtB - \vte|^k (\vtB - \vte) }{\vte} \D \sigma_x }.
				\label{a4} 
			\end{align}
			for any $\varphi \in C^1(S_T \times \Ov{\Omega})$, $\varphi \geq 0$;
			
			\item {\bf Energy balance:}
			\begin{align}  
				- &\intST{ \partial_t \psi	\intO{ \left[ \frac{1}{2} \vre |\vue|^2 + \vre e  \right] } } + \frac{1}{\ep} \intST{ \psi  \int_{\partial \Omega} |\vte - \vtB|^k (\vte - \vtB)  \D \sigma_x
				 }  \br
				&= 
				\intST{ \psi \intO{ \vre \vue \cdot \Grad G  } }
				\label{a5}
			\end{align}
			for any $\psi \in C^1(S_T)$
		\end{itemize}
cf. \cite[Section 2.2]{FeMuNoPo}.	

\subsection{Existence of approximate solutions}	

Our aim is use the existence result proved in \cite[Theorem 1]{FeMuNoPo} to obtain the approximate solutions 
$(\vre, \vte, \vue)_{\ep > 0}$. To perform this step some comments are in order. In comparison with 
\cite{FeMuNoPo}, the present problem features the following new ingredients: 

\begin{itemize}
	
	\item The action of the Coriolis force in the momentum equation \eqref{a3}.
	\item The function $\vtB$ in \eqref{A1} is time dependent whereas its counterpart $\Theta_0$ in 
	\cite{FeMuNoPo} depends only on $x$.
	\item The exponent $k$ in \eqref{A1} equals zero in \cite{FeMuNoPo}.

	\end{itemize}

It is easy to check that the existence proof in \cite{FeMuNoPo} can be modified to accommodate the above changes as soon as suitable {\it a priori} bounds similar to those in \cite[Section 2.4]{FeMuNoPo} are established. 
To see this, we start with the energy balance \eqref{a5} yielding 
\begin{align}
	\frac{1}{\ep} &\intST{ \int_{\partial \Omega} |\vte - \vtB|^k (\vte - \vtB)  \D \sigma_x
	} = \intST{ \intO{ \vre \vue \cdot \Grad G } } = - \intST{ \intO{ \vre \partial_t G } } \br 
	&\leq M_0 \| \partial_t G \|_{L^\infty(S_T \times \Omega)} ,\ 
	\mbox{where}\ M_0 = \intO{ \vr }. 
	\label{A9}
\end{align}
As $\vte > 0$ a.a., \eqref{A9} yields the bound 
\begin{equation} \label{A10}
	\| \vte \|_{L^{k+1}(S_T \times \partial \Omega)} \aleq 1
\end{equation}
in terms of the data and uniform for $\ep \to 0$. Consequently, the entropy inequality \eqref{a4} gives rise
to the bound on the entropy production rate 
\begin{equation} \label{A10a}
\intST{ \intO{ \frac{1}{\vte} \left[ \mathbb{S} : \Ds \vue - 
				\frac{\vc{q} \cdot \Grad \vte }{\vte} \right] } } \leq c(\ep, {\rm data})
			\end{equation} 
and the remaining estimates are obtained exactly as in \cite[Section 2.4]{FeMuNoPo}. Note that the right--hand side of \eqref{A10a} may blow up for $\ep \to 0$. 

With the necessary {\it a priori} bounds at hand, we obtain a family of approximate solutions $(\vre, \vte, \vue)_{\ep > 0}$ exactly as in \cite{FeMuNoPo}. 

\begin{Proposition} [{\bf Approximate solutions}] \label{AP1}
	
	In addition to the hypotheses of Theorem \ref{MT1}, let 
	\begin{equation} \label{A11}
	\ep > 0,\ 6 < k + 1 = \beta, \ M_0 > 0
	\end{equation}
	be given.
	
	Then the approximate problem \eqref{a1} -- \eqref{a5} admits a solution $(\vre, \vte, \vue)$.

	\end{Proposition}

\subsection{Approximate ballistic energy balance}

Let $\tvt \in C^1(S_T \times \Ov{\Omega})$ satisfy \eqref{w9}.
Choosing $\varphi(t,x) = \psi(t) \tvt(t,x)$, where $\psi \in C^1(S_T)$, $\psi \geq 0$, 
as a test function in the approximate entropy inequality \eqref{a4} and adding the resulting integral to the energy balance \eqref{a5}, we deduce
\begin{align}  
	- &\intST{ \partial_t \psi	\intO{ \left[ \frac{1}{2} \vre |\vue|^2 + \vre e - \tvt \vre s \right] } } 
	+ \intST{ \psi
		\intO{ \frac{\tvt}{\vte}	 \left[ \mathbb{S}: \Ds \vue - \frac{\vc{q} \cdot \Grad \vte }{\vte} \right] } }  \br
	&+ \frac{1}{\ep} \intST{ \psi \int_{\partial \Omega} \frac{ |\vte - \vtB|^{k + 2} }{\vte} \ \D \sigma_x }
	\br
	&\leq 
	\intST{ \psi \intO{ \left[ \vre \vue \cdot \Grad G - \vre s \vue \cdot \Grad \tvt - \frac{\vc{q}}{\vte} \cdot \Grad \tvt 
	- \partial_t \tvt \vre s 
	\right] } }.
	\label{a6}
\end{align}
Inequality \eqref{a6} is obviously a counter part of the ballistic energy balance \eqref{w8} and will be used in the forthcoming part to deduce the necessary bounds on the family of approximate solutions.

\section{Uniform bounds}
\label{U}

In order to perform the limit $\ep \to 0$ in the family of approximate solutions obtained in Proposition \ref{AP1}, we need uniform bounds independent of $\ep$.

\subsection{Mass conservation}

Obviously, as the total mass of the fluid is conserved, we get 
\begin{equation} \label{U1}
	M_0 = \intO{ \vre(t,\cdot) } \ \mbox{for all}\ t \in S_T \ 
	\Rightarrow \ \sup_{t \in S_T} \| \vre (t, \cdot)\|_{L^1(\Omega)} \aleq 1.
	\end{equation}

\subsection{Energy estimates} 
\label{EE}

As both $\partial \Omega$ and the boundary data $\vtB$ are smooth, we may suppose that 
\begin{equation} \label{U2}
	\Del \vtB(t, \cdot) = 0 \ \mbox{in}\ \Omega \ \mbox{for any}\ t \in S_T.
	\end{equation}
Choosing $\psi = 1$, $\tvt = \vtB$ in the ballistic energy inequality \eqref{a6} we get 
\begin{align}
\intST{ 
	\intO{ \frac{\vtB}{\vte}	 \left[ \mathbb{S}(\Ds \vue) : \Ds \vue + \frac{\kappa (\vte) |\Grad \vte|^2 }{\vte} \right] } }	 \br
+\frac{1}{\ep} \intST{  \int_{\partial \Omega} \frac{ |\vte - \vtB|^{k + 2} }{\vte} \ \D \sigma_x } \br 
\leq \intST{  \intO{ \left[ \vre \vue \cdot \Grad G - \vre s(\vre, \vte) \vue \cdot \Grad \vtB + \frac{\kappa (\vte) \Grad \vte }{\vte} \cdot \Grad \vtB - \partial_t \tvt \vre s(\vre, \vte) \right] } }.
\label{U3}	
	\end{align}
By virtue of hypothesis \eqref{w16} and Korn's inequality, we obtain 
\[
\| \vue \|_{W^{1,2}_0(\Omega; R^d)}^2 \aleq \intO{ \frac{\vtB}{\vte} \mathbb{S}(\Ds \vue) : \Ds \vue }. 
\]
Moreover, again by virtue of \eqref{w16}, 
\[
\intO{ \left[ | \Grad \vte^{\frac{\beta}{2}} |^2 + |\Grad \log (\vte) |^2 \right] } 
\aleq \intO{ \frac{\vtB}{\vte}	\frac{\kappa (\vte) |\Grad \vte|^2 }{\vte}  }.
\]
By  Poincar\` e inequality (see e.g. Theorem 4.4.6 in \cite{Ziemer}) we obtain that
\[
\intO{  |  \vte^{\frac{\beta}{2}} |^2  } 
\aleq \int_{\partial \Omega}|\vte^{\frac{\beta}{2}}|^2\ \D \sigma_x+\intO{\left|\vte^{\frac{\beta}{2}}- \int_{\partial \Omega} \vte^{\frac{\beta}{2}}\ \D \sigma_x\right|^2}
\aleq \int_{\partial \Omega}|\vte^{\frac{\beta}{2}}|^2\ \D \sigma_x+\intO{ | \Grad \vte^{\frac{\beta}{2}} |^2}
\]
as well as
\[
\intO{  | \log (\vte) |^2  } 
\aleq \int_{\partial \Omega}|\log (\vte)|^2\ \D \sigma_x+\intO{ | \Grad \log (\vte) |^2}.
\]
Collecting the last three inequalities, hypothesis \eqref{A11} and estimating the boundary terms
\[
|\vte^{\frac{\beta}{2}}|^2+|\log (\vte)|^2\aleq \frac{ |\vte - \vtB|^{k + 2} }{\vte}
\]
gives 
\[
\left\| \vte^{\frac{\beta}{2}} \right\|_{W^{1,2}(\Omega)}^2  + 
\left\| \log(\vte) \right\|_{W^{1,2}(\Omega)}^2 
\aleq \left( 1 + \frac{1}{2\ep}   \int_{\partial \Omega} \frac{ |\vte - \vtB|^{k + 2} }{\vte} \ \D \sigma_x + 
\intO{ \frac{\vtB}{\vte}	\frac{\kappa (\vte) |\Grad \vte|^2 }{\vte}  } \right).
\]
Note that as $\ep$ is small, thus the coefficient at the boundary term may appear without violating the estimates, but aiming to match the formulation of approximate ballistic energy balance. 
Gathering the previous observations, we may infer that 
\begin{align}
\intST{ \left[ \| \vue \|_{W^{1,2}_0(\Omega; R^d)}^2 + 	\left\| \vte^{\frac{\beta}{2}} \right\|_{W^{1,2}(\Omega)}^2  + 
	\left\| \log(\vte) \right\|_{W^{1,2}(\Omega)}^2 \right] } \br 
+\frac{1}{\ep} \intST{  \int_{\partial \Omega} \frac{ |\vte - \vtB|^{\beta + 1} }{\vte} \ \D \sigma_x } \br 
\aleq \left(1 + \left| \intST{  \intO{ \left[ \vre \vue \cdot \Grad G - \vre s(\vre, \vte) \vue \cdot \Grad \vtB + \frac{\kappa (\vte) \Grad \vte }{\vte} \cdot \Grad \vtB - \partial_t \tvt \vre s(\vre, \vte)\right] } } \right| \right).
	\label{U4}
	\end{align}

Now, as $\vre$, $\vue$ solve the equation of continuity \eqref{w3}, 
\[
\intST{ \intO{ \vre \vue \cdot \Grad G }} = - \intST{ \vre \partial_t G } \leq c(M_0, G).
\]
In addition, denoting 
\[
\mathcal{K}(\vt) = \int_1^\vt \frac{\kappa (z) }{z} \ \D z, 
\]
we obtain, by virtue of \eqref{U2}, 
\[
\intO{ \frac{\kappa (\vte) \Grad \vte }{\vte} \cdot \Grad \vtB } = 
\intO{ \Grad \mathcal{K} (\vte) \cdot \Grad \vtB } = \int_{\partial \Omega} \mathcal{K} (\vte) \Grad \vtB \cdot 
\vc{n} \D \sigma_x.
\] 
Consequently, as $\kappa$ satisfies hypothesis \eqref{w16}, we conclude 
\[
\left| \intO{ \frac{\kappa (\vte) \Grad \vte }{\vte} \cdot \Grad \vtB } \right| 
= \left| \int_{\partial \Omega} \mathcal{K} (\vte) \Grad \vtB \cdot 
\vc{n} \D \sigma_x \right| \aleq \left(1+ \intST{  \int_{\partial \Omega} \frac{ |\vte - \vtB|^{\beta + 1} }{\vte} \ \D \sigma_x}\right).
\]
Thus inequality \eqref{U4} reduces to 
\begin{align}
	\intST{ \left[ \| \vue \|_{W^{1,2}_0(\Omega; R^d)}^2 + 	\left\| \vte^{\frac{\beta}{2}} \right\|_{W^{1,2}(\Omega)}^2  + 
		\left\| \log(\vte) \right\|_{W^{1,2}(\Omega)}^2 \right] } \br 
	+\frac{1}{\ep} \intST{  \int_{\partial \Omega} \frac{ |\vte - \vtB|^{\beta + 1} }{\vte} \ \D \sigma_x } \br 
	\aleq \left(1 + \intST{  \intO{ \left(| \vre s(\vre, \vte) \vue \cdot \Grad \vtB | + |\partial_t \tvt \vre s(\vre, \vte)|\right) } }  \right).
	\label{U5}
\end{align}
In accordance with hypothesis \eqref{w12}, we decompose the entropy as 
\[
\vre s(\vre, \vte)  = \vre \mathcal{S} \left( \frac{\vre}{\vte^{\frac{3}{2}} } \right) + \frac{4a}{3} \vte^3.
\]
Consequently, the radiation component may be handled as 
\[
\intO{ | \vte^3 \vue \cdot \Grad \vtB | } \leq \delta \| \vue \|_{L^2(\Omega; R^d)}^2 + 
c(\delta, \vtB) \intO{ \vte^6 } 
\]
for any $\delta > 0$. Consequently, as $\beta > 6$, this term can be absorbed by the left--hand side of 
\eqref{U5} yielding 
\begin{align}
	\intST{ \left[ \| \vue \|_{W^{1,2}_0(\Omega; R^d)}^2 + 	\left\| \vte^{\frac{\beta}{2}} \right\|_{W^{1,2}(\Omega)}^2  + 
		\left\| \log(\vte) \right\|_{W^{1,2}(\Omega)}^2 \right] } \br 
	+\frac{1}{\ep} \intST{  \int_{\partial \Omega} \frac{ |\vte - \vtB|^{\beta + 1} }{\vte} \ \D \sigma_x }  
	\aleq \left(1 + \intST{  \intO{ \left| \vre \mathcal{S} \left( \frac{\vre}{\vte^{\frac{3}{2}} } \right) \vue \cdot \Grad \vtB \right|+\left|\vre \mathcal{S} \left( \frac{\vre}{\vte^{\frac{3}{2}} } \right) \partial_t \tvt \right|  } }  \right).
	\label{U6}
\end{align}

Finally, following the arguments of \cite[Section 4.4]{FeiSwGw}, we make use of the Third law of thermodynamics 
enforced through hypothesis \eqref{w14}. Specifically, if 
\[
\frac{\vr}{\vt^{\frac{3}{2}}} < r \ \mbox{meaning}\ \vr < r \vt^{\frac{3}{2}}, 
\]
we get, by virtue of \eqref{w15}, 
\begin{equation} \label{U7}
	0 \leq \vr \mathcal{S} \left( \frac{\vr}{\vt^{\frac{3}{2}}} \right) 
	\aleq \left(1 + r \vt^{\frac{3}{2}} \left[ \log^+ (r \vt^{\frac{3}{2}} ) + \log^+(\vt) \right] \right). 
\end{equation}
Consequently, we deduce from \eqref{U6}, 
\begin{align}
	\intST{ \left[ \| \vue \|_{W^{1,2}_0(\Omega; R^d)}^2 + 	\left\| \vte^{\frac{\beta}{2}} \right\|_{W^{1,2}(\Omega)}^2  + 
		\left\| \log(\vte) \right\|_{W^{1,2}(\Omega)}^2 \right] } 
	+\frac{1}{\ep} \intST{  \int_{\partial \Omega} \frac{ |\vte - \vtB|^{\beta + 1} }{\vte} \ \D \sigma_x }  \br
	\aleq \left(\Lambda(r) + \intST{ \intO{ \mathds{1}_{ \left\{ \frac{\vre}{\vte^{\frac{3}{2}}} \geq r \right\} }  \left| \vre \mathcal{S} \left( \frac{\vre}{\vte^{\frac{3}{2}} } \right) \vue \cdot \Grad \vtB \right|  } }  
	 + \intST{ \intO{ \mathds{1}_{ \left\{ \frac{\vre}{\vte^{\frac{3}{2}}} \geq r \right\} }  \left| \vre \mathcal{S} \left( \frac{\vre}{\vte^{\frac{3}{2}} } \right)\partial_t \tvt \right|  } }  
	\right),
	\label{U8}
\end{align}
where $\Lambda(r) \to \infty$ as $r \to \infty$. 

Now, again by hypothesis \eqref{w14}, 
\[
0 \leq \mathds{1}_{ \left\{ \frac{\vre}{\vte^{\frac{3}{2}}} \geq r \right\} }\mathcal{S} \left( \frac{\vre}{\vte^{\frac{3}{2}} } \right) \leq \mathcal{S}(r) \to 0 \ \mbox{as}\ r \to 
\infty.
\]
In an analogue way we treat the term $\intST{ \intO{  \partial_t \tvt \vre s  } } $.
Going back to \eqref{U8} we conclude 
\begin{align}
	\intST{ \left[ \| \vue \|_{W^{1,2}_0(\Omega; R^d)}^2 + 	\left\| \vte^{\frac{\beta}{2}} \right\|_{W^{1,2}(\Omega)}^2  + 
		\left\| \log(\vte) \right\|_{W^{1,2}(\Omega)}^2 \right] } \br 
	+\frac{1}{\ep} \intST{  \int_{\partial \Omega} \frac{ |\vte - \vtB|^{\beta + 1} }{\vte} \ \D \sigma_x }  
	\aleq \left(\Lambda(r) + \mathcal{S}(r) \intST{ \intO{  \left( \left| \vre \vue \right|  +\vre\right)} }  \right), \br 
	\Lambda(r) \to \infty,\ \mathcal{S}(r) \to 0 \ \mbox{as}\ r \to \infty.
	\label{U9}
\end{align}

\subsection{Pressure estimates}
\label{PE}

To close the estimates we have to control the density in terms of the integrals on the right--hand side of 
\eqref{U9}. To this end, we use the nowadays standard pressure estimates obtained via Bogovskii operator. 
Specifically, we use the quantity
\[
\varphi (t,x) = \mathcal{B} \left[ \vre^\omega - \frac{1}{|\Omega|} \intO{ \vre^\omega } \right], \omega > 0,
\]
as a test function in the momentum equation \eqref{a3}. Here $\mathcal{B}$ denotes the operator enjoying 
the following properties: 
\begin{itemize}
	\item 
	\[
	\mathcal{B} : L^q_0 (\Omega) \equiv \left\{ v \in L^q(\Omega) \ \Big|\ \intO{ v } = 0 \right\} \to 
	W^{1,q}_0 (\Omega; R^d) ,\ 1 < q < \infty;
	\]
	\item
	\[
	\Div \mathcal{B} [v] = v; 
	\]
	\item
	if $v = \Div \vc{g}$, with $\vc{g} \in L^q(\Omega; R^d)$, $\Div \vc{g} \in L^r(\Omega)$, 
	$\vc{g} \cdot \vc{n}|_{\partial \Omega} = 0$, then 
	\[
	\| \mathcal{B}[ \Div \vc{g} ]\|_{L^r(\Omega; R^d)} \aleq \| \vc{g} \|_{L^r(\Omega; R^d)},
	\]
	
	\end{itemize}
see Galdi \cite[Chapter 3]{GAL} or Gei\ss ert, Heck, and Hieber \cite{GEHEHI}. After a straightforward 
manipulation (see e.g. \cite{FMPS}), we obtain 
\begin{align}
	\int_{S_T} &\intO{ p(\vre, \vte) \vre^\omega } \dt = \intST{ \frac{1}{|\Omega|} \left( \intO{ \vre^\omega }  \right) \left( \intO{ p(\vre, \vte)  } \right) } \br 
	&-\intST{ \intO{ \vre (\vue \otimes \vue): \Grad \mathcal{B} \left[ \vre^\omega - \frac{1}{|\Omega|} \intO{\vre^\omega } \right] } 
} \br 
&+ \intST{\intO{\vre (\bfomega \times \vue) \cdot  \mathcal{B} \left[ \vre^\omega - \frac{1}{|\Omega|} \intO{ \vre^\omega } \right]      }}
\br
	&+ \intST{ \intO{ \mathbb{S}(\vte, \Ds \vue) : \Grad \mathcal{B} \left[ \vre^\omega - \frac{1}{|\Omega|} \intO{ \vre^\omega } \right] } 
} \br & - \intST{ \intO{ \vre \Grad G \cdot \mathcal{B} \left[ \vre^\omega - \frac{1}{|\Omega|} \intO{ \vre^\omega }  \right] } } \br 
&+ \intST{ \intO{	\vre \vue \cdot \mathcal{B}[\Div ( \vre^\omega \vue) ] } } \br 
&+ (\omega - 1) \intST{ \intO{ \vre \vue \cdot \mathcal{B} \left[ \vre^\omega \Div \vue - 
	\frac{1}{|\Omega|} \intO{ \vre^\omega \Div \vue } \right] } }.
	\label{U10}	
\end{align}
Since the total mass $M_0$ is constant, the smoothing properties of $\mathcal{B}$ yield 
\[
\left\| \mathcal{B} \left[ \vre^\omega - \frac{1}{|\Omega|} \intO{ \vre^\omega } \right] \right\|_{L^\infty(S_T 
	\times \Omega; R^d) } \leq c(M_0) \ \mbox{as soon as}\ \omega < \frac{1}{d}.
\]
Moreover, in accordance with hypotheses \eqref{M2}--\eqref{w11}, 
\[
\vr^{\frac{5}{3}} + \vt^4 \aleq p(\vr, \vt) \aleq \vr^{\frac{5}{3}} + \vt^4 + 1.
\]
In view of these facts, inequality \eqref{U10} gives rise to 
\begin{align}
	\int_{S_T} &\intO{ \vre^{\frac{5}{3} + \omega} } \dt \leq c(M_0) \Big( 1 +  \intST{ \intO{ \vte^4 } } \br 
	&-\intST{ \intO{ \vre (\vue \otimes \vue): \Grad \mathcal{B} \left[ \vre^\omega - \frac{1}{|\Omega|} \intO{\vre^\omega } \right] } 
	} \br 
	&+ \intST{\intO{\vre (\bfomega \times \vue) \cdot  \mathcal{B} \left[ \vre^\omega - \frac{1}{|\Omega|} \intO{ \vre^\omega } \right]      }}
	\br
	&+ \intST{ \intO{ \mathbb{S}(\vte, \Ds \vue) : \Grad \mathcal{B} \left[ \vre^\omega - \frac{1}{|\Omega|} \intO{ \vre^\omega } \right] } 
	} \br  
	&+ \intST{ \intO{	\vre \vue \cdot \mathcal{B}[\Div ( \vre^\omega \vue) ] } } \br 
	&+ (\omega - 1) \intST{ \intO{ \vre \vue \cdot \mathcal{B} \left[ \vre^\omega \Div \vue - 
			\frac{1}{|\Omega|} \intO{ \vre^\omega \Div \vue } \right] } } \Big).
	\label{U11}	
\end{align}

The following steps will be performed for $d=3$. Obviously even better estimates can be obtained if $d=2$. 
First, 
\begin{align}
	&\left| \intST{ \intO{ \vre (\vue \otimes \vue): \Grad \mathcal{B} \left[ \vre^\omega - \frac{1}{|\Omega|} \intO{ \vre^\omega } \right] } } \right| \br
	&\quad \aleq \intST{ \| \vre \|_{L^\gamma(\Omega)} \| \vue \|^2_{L^6(\Omega; R^3)} \| \vre^\omega \|_{L^q(\Omega)}  
}	\br &\quad \aleq \sup_{t \in S_T} \| \vre \|_{L^\gamma (\Omega)} \intST{ \| \vue \|^2_{W^{1,2}(\Omega; R^3)} } \sup_{t \in S_T} \| \vr^\omega \|_{L^q (\Omega)}
	\dt,
	\nonumber
\end{align}
where 
\[
q = \frac{3 \gamma}{2 \gamma - 3} > 1 \ \mbox{provided}\ \gamma > \frac{3}{2}. 
\]
Fixing 
\begin{equation} \label{U12}
\gamma = \frac{5}{3},\ \omega = \frac{3 \gamma}{2 \gamma - 3} = \frac{1}{15}
\end{equation}
and using the fact that the total mass $M_0$ is conserved, we get 
\begin{align}
	&\left| \intST{ \intO{ \vre (\vue \otimes \vue): \Grad \mathcal{B} \left[ \vre^\omega - \frac{1}{|\Omega|} \intO{ \vre^\omega } \right] } } \right| 	\br &\quad \leq c(M_0) \sup_{t \in S_T} \| \vre \|_{L^\gamma (\Omega)} \intST{ \| \vue \|^2_{W^{1,2}(\Omega; R^3)} }.
	\nonumber
\end{align}
Seeing that the integral containing the Coriolis force can be controlled in a similar way we may rewrite \eqref{U11} in the form
\begin{align}
	\int_{S_T} &\intO{ \vre^{\frac{5}{3} + \omega} } \dt \leq c(M_0) \Big( 1 +  \intST{ \intO{ \vte^4 } } \br 
	&+ \sup_{t \in S_T} \| \vre \|_{L^\gamma (\Omega)} \intST{ \| \vue \|^2_{W^{1,2}(\Omega; R^3)} } \br
	&+ \intST{ \intO{ \mathbb{S}(\vte, \Ds \vue) : \Grad \mathcal{B} \left[ \vre^\omega - \frac{1}{|\Omega|} \intO{ \vre^\omega } \right] } 
	} \br  
	&+ \intST{ \intO{	\vre \vue \cdot \mathcal{B}[\Div ( \vre^\omega \vue) ] } } \br 
	&+ (\omega - 1) \intST{ \intO{ \vre \vue \cdot \mathcal{B} \left[ \vre^\omega \Div \vue - 
			\frac{1}{|\Omega|} \intO{ \vre^\omega \Div \vue } \right] } } \Big).
	\label{U13}	
\end{align}
In a similar way, we get
\begin{align}
	\left| \intST{ \intO{ \vre \vue \cdot \mathcal{B}[\Div (\vre^\omega \vue)] } } \right|	 \aleq 
	\intST{ \| \vre \|_{L^\gamma (\Omega)} \| \vue \|_{L^6(\Omega; R^3)} \| \vre^\omega \vue \|_{L^q(\Omega; R^3)} },
	\nonumber
\end{align}
where
\[
\frac{1}{\gamma}+ \frac{1}{6} + \frac{1}{q} = 1.
\]
In addition, 
\[
\| \vre^\omega \vue \|_{L^q(\Omega; R^3)} \leq \| \vue \|_{L^6(\Omega; R^3)} \| \vre^\omega \|_{L^p(\Omega)}, 
\mbox{where}\ 
\frac{1}{q} = \frac{1}{6} + \frac{1}{p};  
\]
whence 
\begin{align}
	\left| \int_{\tau}^{\tau + 1} \intO{ \vr \vu \cdot \mathcal{B}[\Div (\vr^\alpha \vu)] } \right|	 \leq 
	c(M) \sup_{t \in (\tau, \tau + 1)} \| \vr \|_{L^{\frac{5}{3}} (\Omega)} \int_{\tau}^{\tau + 1} \| \vu \|^2_{W^{1,2}(\Omega; R^3)} \dt	
	\nonumber
\end{align}
as soon as  \eqref{U12} holds.

Finally, 
\begin{align}
	&\left|  \intST{ \intO{ \vre \vue \cdot \mathcal{B} \left[ \vre^\omega \Div \vue - 
		\frac{1}{|\Omega|} \intO{ \vre^\omega \Div \vue } \right] } } \right| \br &\quad \leq 
	\intST{ \| \vre \|_{L^\gamma (\Omega)} \| \vue \|_{L^6(\Omega; R^3)} \left\| \mathcal{B} \left[ \vre^\omega \Div \vue - 
	\frac{1}{|\Omega|} \intO{ \vre^\omega \Div \vue } \right] \right\|_{L^q(\Omega; R^3)} },
	\nonumber
\end{align}
where
\[
\frac{1}{\gamma} + \frac{1}{6} + \frac{1}{q} = 1.
\]
Here, 
\[
\left\| \mathcal{B} \left[ \vre^\omega \Div \vue - 
\frac{1}{|\Omega|} \intO{ \vre^\omega \Div \vue } \right] \right\|_{L^q(\Omega; R^3)} \aleq 
\| \vre^\omega \Div \vue \|_{L^r(\Omega; R^3)},\ q = \frac{3 r}{3 - r}, 
\]
and 
\[
\| \vre^\omega \Div \vue \|_{L^r(\Omega; R^3)} \leq \| \vue \|_{W^{1,2}(\Omega; R^3)} \| \vre^\omega \|_{L^p(\Omega)},\ 
\mbox{with}\ \frac{1}{2} + \frac{1}{p} = \frac{1}{r}.
\]
Thus using \eqref{U12} we may infer that 
\begin{align}
	&\left|  \intST{ \intO{ \vre \vue \cdot \mathcal{B} \left[ \vre^\omega \Div \vue - 
		\frac{1}{|\Omega|} \intO{ \vre^\omega \Div \vue } \right] } } \right| \br &\quad \leq 
	c(M_0) \sup_{t \in S_T} \| \vre \|_{L^{\frac{5}{3}} (\Omega)} \intST{ \| \vue \|^2_{W^{1,2}(\Omega; R^3)} }.
	\nonumber
\end{align}
Going back to \eqref{U13} and summarizing the previous estimates we conclude 
\begin{align}
	\int_{S_T} &\intO{ \vre^{\frac{5}{3} + \omega} } \dt \leq c(M_0) \left( 1 + \intST{ \intO{ \vte^4 }  } \right. \br 
	&+	\sup_{t \in S_T} \| \vre \|_{L^{\frac{5}{3}} (\Omega)} \intST{ \| \vue \|^2_{W^{1,2}(\Omega; R^3)}
}
	\br 
	&+ \left. \intST{ \intO{ \mathbb{S}(\vte, \Ds \vue) : \Grad \mathcal{B} \left[ \vre^\omega - \frac{1}{|\Omega|} \intO{ \vre^\omega } \right] } 
} \right),\ \mbox{where}\  \omega = \frac{1}{15}.
	\label{U14}	
\end{align}	
The last step is estimating 
\begin{align} 
	&\intO{ \mathbb{S} (\vte, \Ds \vue) : \Grad \mathcal{B} \left[ \vre^\omega - \frac{1}{|\Omega|} \intO{ \vr^\omega } \right] } \br 
	&\quad \leq  (1 + \| \vte \|_{L^4(\Omega)} ) \| \vue \|_{W^{1,2}(\Omega; R^3)} \left\| 
	\Grad \mathcal{B} \left[ \vre^\omega - \frac{1}{|\Omega|} \intO{ \vre^\omega } \right] \right\|_{L^4(\Omega; R^3)} \br 
	&\quad \leq c(M) (1 + \| \vte \|_{L^4(\Omega)} ) \| \vue \|_{W^{1,2}(\Omega; R^3)}.
	\nonumber 
\end{align} 
We therefore conclude the pressure estimates:
\begin{align}
	\int_{S_T} &\intO{ \vre^{\frac{5}{3} + \omega} } \dt \leq c(M_0) \left[ 1 + \intST{ \intO{ \vte^4 }  } \right. \br 
	&+	\left. \left(1 + \sup_{t \in S_T} \| \vre \|_{L^{\frac{5}{3}} (\Omega)} \right) \intST{ \| \vue \|^2_{W^{1,2}(\Omega; R^3)}
}   \right],\ \omega = \frac{1}{15}. 
	\label{U15}	
\end{align}	

\subsection{Uniform bounds for $\ep \to 0$}

As $\beta > 6$, we deduce from the inequalities \eqref{U9}, \eqref{U15} that 
\begin{align}
	\intST{ \intO{ \vte^4 }} &\aleq \left(1 +  \intST{ \| \vte^{\frac{\beta}{2}} \|^2_{W^{1,2}(\Omega)} 
} \right) \br &\aleq  \left(1 + \intST{ \intO{ \vre |\vue| } } \right)
\nonumber	
\end{align}
provided we fix $r = 1$ in \eqref{U9}. Furthermore, 
\begin{align}
\intST{ \intO{ \vre |\vue| } } \leq 
\frac{1}{2} \intST{ \intO{ \vre } } + \frac{1}{2} \intST{ \intO{ \vre |\vue|^2 }} \br 
\leq \frac{1}{2}T M_0 + \frac{1}{2} \sup_{t \in S_T} \| \vre \|_{L^{\frac{5}{3}}(\Omega)} \intST{ \| \vue \|^2_{L^5(\Omega; R^3)}} \br 
\leq c(M_0) \left(1 +  \sup_{t \in S_T} \| \vre \|_{L^{\frac{5}{3}}(\Omega)} \intST{ \| \vue \|^2_{W^{1,2}(\Omega; R^3)}} \right)
			\nonumber
\end{align}
Consequently, inequality \eqref{U15} reduces to 
\begin{align}
	\int_{S_T} &\intO{ \vre^{\frac{5}{3} + \omega} } \dt \leq c(M_0) \left[ 1 +  \left(1 + \sup_{t \in S_T} \| \vre \|_{L^{\frac{5}{3}} (\Omega)} \right) \intST{ \| \vue \|^2_{W^{1,2}(\Omega; R^3)}
	}   \right],\ \omega = \frac{1}{15}. 
	\label{U16}	
\end{align}	
Next, going back to \eqref{U9} we get 
\[
\intST{ \| \vue \|^2_{W^{1,2}(\Omega; R^3)} }  \aleq \left( \mathcal{S}(r) \intST{ \intO{ (\vre |\vue| +\vre)} }
+ \Lambda (r) \right)
\]
where, by means of the standard Sobolev embedding theorem, 
\[
\intO{ \vre |\vue| } \leq  \| \sqrt{\vre} \|_{L^2(\Omega)} \| \sqrt{\vre} \|_{L^3(\Omega)} \| \vue \|_{L^6(\Omega; R^3)} 
\leq c(M_0) \| \sqrt{\vre} \|_{L^3(\Omega)} \| \vue \|_{W^{1,2}(\Omega; R^3)}.
\]
Consequently, 
\begin{equation} \label{U17}
\intST{ \| \vue \|^2_{W^{1,2}(\Omega; R^3)} }  \aleq \left( \mathcal{S}(r) \intST{ \| \vre \|_{L^{\frac{3}{2}}(\Omega) } }
+ \Lambda (r) \right)
\end{equation}
Now, introducing the total energy of the system, 
\[
E(\vr, \vt, \vu) = \frac{1}{2} \vr |\vu|^2 + \vr e(\vr, \vt) 
\]
we first observe that 
\begin{equation} \label{U18}
\sup_{t \in S_T} \intO{ E(\vre, \vte, \vue) } \aleq \left( 1 + \intST{ \intO{  E(\vre, \vte, \vue) } } \right).
	\end{equation}
The estimate \eqref{U18} follows from the mean value theorem and the ballistic energy inequality \eqref{a6}. 
Indeed, in view of the uniform bounds established in Section \ref{EE}, we first deduce \eqref{U18} 
for the ballistic energy 
\[
E(\vre, \vte, \vue) - \vtB \vre s(\vre, \vte),
\] 
and then use \eqref{w15} to observe that the entropy part $\vt_B \vre s(\vre, \vte)$ is a lower order perturbation.

Now, we estimate the kinetic energy using \eqref{U17}, 
\begin{align}
	\int_{S_T} & \intO{ \vre |\vue|^2 } \dt \leq \sup_{t \in S_T} \| \vre \|_{L^{\frac{3}{2}}(\Omega)} 
	\intST{ \| \vue \|^2_{L^6(\Omega; R^3)} } \br &\leq c \sup_{t \in S_T} \| \vre \|_{L^{\frac{3}{2}}(\Omega)} 
	\int_{\tau}^{\tau + 1} \| \vue \|^2_{W^{1,2}(\Omega; R^3)} \dt \br 
	&\aleq \Lambda (r) \sup_{t \in S_T} \| \vre \|_{L^{\frac{3}{2}}(\Omega)} + 
	\mathcal{S}(r) \sup_{t \in S_T} \| \vre \|_{L^{\frac{3}{2}}(\Omega)}\intST{ \| \vre \|_{L^{\frac{3}{2}}(\Omega)} }. \nonumber
\end{align}
In addition, by interpolation, 
\begin{equation} \label{U19}
\| \vre \|_{L^{\frac{3}{2}}(\Omega)} \leq \| \vre \|_{L^{\frac{5}{3}}(\Omega)}^{\frac{5}{6}} \| \vre \|_{L^1(\Omega)}^{\frac{1}{6}}.
\end{equation}
Consequently, 
\begin{align}
	\int_{S_T} & \intO{ \vre |\vue|^2 } \dt \leq \br 
	&\aleq c(M_0) \Lambda (r) \sup_{t \in S_T} \| \vre \|_{L^{\frac{5}{3}}(\Omega)}^{\frac{5}{6}} + 
	\mathcal{S}(r) \sup_{t \in S_T} \| \vre \|_{L^{\frac{5}{3}}(\Omega)}^{\frac{5}{6}} \intST{ \| \vre \|_{L^{\frac{5}{3}}(\Omega)}^{\frac{5}{6}} }. \label{U20}
\end{align}

Combining \eqref{U16}, \eqref{U17}, \eqref{U19} we get 
\begin{align}
	\int_{S_T} &\intO{ \vre^{\frac{5}{3} + \omega} } \dt \leq c(M_0) \left[ 1 +  \left(1 + \sup_{t \in S_T} \| \vre \|_{L^{\frac{5}{3}} (\Omega)} \right) \intST{ \| \vue \|^2_{W^{1,2}(\Omega; R^3)}
	}   \right]\br 
&\leq c(M_0) \left[   1 +  \left(1 + \sup_{t \in S_T} \| \vre \|_{L^{\frac{5}{3}} (\Omega)} \right)\left( \mathcal{S}(r) \intST{ \| \vre \|_{L^{\frac{3}{2}}(\Omega) } }
+ \Lambda (r) \right)  \right] \br 
&\leq c(M_0) \left[   1 +  \left(1 + \sup_{t \in S_T} \| \vre \|_{L^{\frac{5}{3}} (\Omega)} \right)\left( \mathcal{S}(r) \intST{ \| \vre \|_{L^{\frac{5}{3}}(\Omega) }^{\frac{5}{6}} }
+ \Lambda (r) \right)  \right] 
	\label{U21}	
\end{align}	

Interpolating $L^1$ and $L^{\frac{5}{3} + \omega}$ and using boundedness of the total mass we have
\begin{equation} \label{U22}
\intST{ \intO{ \vre^{\frac{5}{3}} } } \leq c(M_0) \left( \intST{ \intO{ \vre^{\frac{5}{3} + \omega } }
} \right)^{\frac{10}{11}} \ \mbox{provided}\ \omega = \frac{1}{15}.
\end{equation}
Thus summing up \eqref{U18}--\eqref{U22} we may infer that 
\begin{align}
	\mbox{\eqref{U18}}\ \Rightarrow \ 
&\sup_{t \in S_T} \intO{ E(\vre, \vte, \vue) } \aleq \left( 1 + \intST{ \intO{  E(\vre, \vte, \vue) } } \right) \br 
&\aleq \left[ 1 + \intST{ \intO{\left(  \| \vue \|_{W^{1,2}_0(\Omega; R^d)}^2 + 	\left\| \vte^{\frac{\beta}{2}} \right\|_{W^{1,2}(\Omega)}^2  + 
		\left\| \log(\vte) \right\|_{W^{1,2}(\Omega)}^2  \right) }} \right. \br 
	&+ \left. \intST{ \intO{ \vre |\vue|^2 }} + \intST{ \intO{ \vre^{\frac{5}{3}} } } \right] \br 
	\mbox{\eqref{U9}}\ \Rightarrow &\aleq  \left[ 1 +  \intST{ \intO{ \vre |\vue|^2 }} + \intST{ \intO{ \vre^{\frac{5}{3}} } } \right]	\br 
\mbox{\eqref{U20}}\ \Rightarrow &\aleq \left[1 +  \intST{ \intO{ \vre^{\frac{5}{3}} } } \right. \br  &+ \left. \Lambda(r) \left( \sup_{t \in S_T} \intO{ E(\vre, \vte, \vue) } \right)^{\frac{1}{2}} 
+ \mathcal{S}(r) \left( \sup_{t \in S_T} \intO{ E(\vre, \vte, \vue) } \right) \right] \br 
\mbox{\eqref{U22}}\ \Rightarrow &\leq c(M_0) \left[1 +  \left( \intST{ \intO{ \vre^{\frac{5}{3} + \omega } }
} \right)^{\frac{10}{11}} \right. \br  &+ \left. \Lambda(r) \left( \sup_{t \in S_T} \intO{ E(\vre, \vte, \vue) } \right)^{\frac{1}{2}} 
+ \mathcal{S}(r) \left( \sup_{t \in S_T} \intO{ E(\vre, \vte, \vue) } \right) \right]  \br 	
\mbox{\eqref{U21}}\ \Rightarrow &\leq c(M_0) \left[\Lambda(r) +  \mathcal{S}^{\frac{10}{11}}(r) \left( \sup_{t \in S_T} \intO{ E(\vre, \vte, \vue) } \right)  \right. \br  &+ \left. \Lambda(r) \left( \sup_{t \in S_T} \intO{ E(\vre, \vte, \vue) } \right)^{\frac{1}{2}} 
+ \mathcal{S}(r) \left( \sup_{t \in S_T} \intO{ E(\vre, \vte, \vue) } \right) \right]. 	
	\label{U23}
	\end{align}
As $\mathcal{S}(r) \to 0$ as $r \to \infty$, we fix $r > 0$ large enough to deduce from \eqref{U23} the desired energy bound 
\begin{equation} \label{U24}
	\sup_{t \in S_T} \intO{ E(\vre, \vte, \vue) } \leq c(M_0).
\end{equation}

\section{Convergence}
\label{C}

Our ultimate goal is to perform the limit in the sequence of approximate solutions $(\vre, \vte, \vue)_{\ep > 0}$ to obtain the existence of the time--periodic solution claimed in Theorem \ref{MT1}. With the energy estimate \eqref{U24} at hand, this is a routine matter nowadays well understood. Indeed the test functions used in the entropy inequality \eqref{w7} are compactly supported thus unaffected 
by the boundary integral in its approximate counterpart \eqref{a4}. Similarly, the approximate ballistic energy \eqref{a6} is in fact stronger than \eqref{w8} due to the penalization 
\begin{equation} \label{C1}
	\frac{1}{\ep} \intST{ \psi \int_{\partial \Omega} \frac{ |\vte - \vtB|^{k+2} }{\vte} \D \sigma_x } \aleq 1,\ \psi \geq 0.
	\end{equation}
In particular, for $\psi = 1$, the above inequality together with the \eqref{a6} yield 
\[
\vte \to \vt \ \mbox{weakly in}\ L^2(0,T; W^{1,2}(\Omega; R^d))
\]
with the limit trace $\vt|_{\partial \Omega} = \vtB$
as required in Theorem \ref{MT1}. 

Consequently, the proof of convergence is exactly the same as in the existence theory elaborated in \cite{ChauFei} with the exception of the strong convergence of the density, the ``initial'' value of which is 
unspecified in the periodic setting. Fortunately, the compactness arguments based on Lions' identity and boundedness of the oscillation defect measure can be modified to accommodate the time periodic
setting exactly as in \cite[Section 9.3]{FeMuNoPo}. Thus the proof of Theorem \ref{MT1} can be completed.

\section{Concluding remarks}

In comparison with \cite{FeMuNoPo}, the available {\it a priori} bounds do not allow to handle a general driving force $\vr \vc{g}$ in the momentum equation. Although the potential case 
$\vc{g} = \Grad G$ is physically relevant, more general (non--potential) forces occur when the fluid is stirred up by the motion of the container. A detailed inspection of the arguments 
in Section \ref{PE} reveals they could be considerably improved in the case $d=2$ due to the Sobolev embedding $W^{1,2} \subset L^q$ for any finite $q$. Similar improvement may also be expected 
in the case the total mass $M_0$ is small, cf. Wang and Wang \cite{WanWan}. We therefore strongly conjecture that the present result can be extended to a general driving force $\vc{g}$ provided 
\begin{itemize}
	\item either $d=2$, 
	\item or 
	\[
	M_0 = \intO{ \vr } 
	\]
	is small enough with respect to the amplitude of $\vc{g}$. 
	
	\end{itemize}
As potentiality of $\vc{g}$ was also used in the estimate \eqref{A9} crucial for boundedness of the approximate sequence, the proof of the above conjecture would require a different kind of approximation scheme.

\def\cprime{$'$} \def\ocirc#1{\ifmmode\setbox0=\hbox{$#1$}\dimen0=\ht0
	\advance\dimen0 by1pt\rlap{\hbox to\wd0{\hss\raise\dimen0
			\hbox{\hskip.2em$\scriptscriptstyle\circ$}\hss}}#1\else {\accent"17 #1}\fi}


\begin{thebibliography}{10}
	
	\bibitem{BEROFO}
	S.~E. Bechtel, F.J. Rooney, and M.G. Forest.
	\newblock Connection between stability, convexity of internal energy, and the
	second law for compressible {N}ewtonian fuids.
	\newblock {\em J. Appl. Mech.}, {\bf 72}:299--300, 2005.
	
	\bibitem{BEL1}
	F.~Belgiorno.
	\newblock Notes on the third law of thermodynamics, {I}.
	\newblock {\em J. Phys. A}, {\bf 36}:8165--8193, 2003.
	
	\bibitem{BEL2}
	F.~Belgiorno.
	\newblock Notes on the third law of thermodynamics, ii.
	\newblock {\em J. Phys. A}, {\bf 36}:8195--8221, 2003.
	
	\bibitem{BreKag2}
	J.~B\v{r}ezina and Y.~Kagei.
	\newblock Decay properties of solutions to the linearized compressible
	{N}avier-{S}tokes equation around time-periodic parallel flow.
	\newblock {\em Math. Models Methods Appl. Sci.}, {\bf 22}(7):1250007, 53, 2012.
	
	\bibitem{BreKag1}
	J.~B\v{r}ezina and Y.~Kagei.
	\newblock Spectral properties of the linearized compressible {N}avier-{S}tokes
	equation around time-periodic parallel flow.
	\newblock {\em J. Differential Equations}, {\bf 255}(6):1132--1195, 2013.
	
	\bibitem{CaiTan}
	H.~Cai and Z.~Tan.
	\newblock Weak time-periodic solutions to the compressible {N}avier-{S}tokes
	equations.
	\newblock {\em Acta Math. Sci. Ser. B (Engl. Ed.)}, {\bf 36}(2):499--513, 2016.
	
	\bibitem{ChauFei}
	N.~Chaudhuri and E.~Feireisl.
	\newblock {N}avier--{S}tokes--{F}ourier system with {D}irichlet boundary
	conditions.
	\newblock {\em {\bf arxiv preprint No. 2106.05315}}, 2021.
	
	\bibitem{DAVI}
	P.~A. Davidson.
	\newblock {\em Turbulence:{A}n introduction for scientists and engineers}.
	\newblock Oxford University Press, Oxford, 2004.
	
	\bibitem{FMPS}
	E.~Feireisl, {\v S}.~Matu{\v s}{\accent23u}-{N}e{\v c}asov{\' a},
	H.~Petzeltov{\' a}, and I.~Stra{\v s}kraba.
	\newblock On the motion of a viscous compressible flow driven by a
	time-periodic external flow.
	\newblock {\em Arch. Rational Mech. Anal.}, {\bf 149}:69--96, 1999.
	
	\bibitem{FeMuNoPo}
	E.~Feireisl, P.~Mucha, A.~Novotn{\' y}, and M.~Pokorn{\' y}.
	\newblock Time periodic solutions to the full {N}avier-{S}tokes-{F}ourier
	system.
	\newblock {\em Arch. Rational. Mech. Anal.}, {\bf 204}:745--786, 2012.
	
	\bibitem{FeNo6A}
	E.~Feireisl and A.~Novotn\'y.
	\newblock {\em Singular limits in thermodynamics of viscous fluids}.
	\newblock Advances in Mathematical Fluid Mechanics. Birkh\"auser/Springer,
	Cham, 2017.
	\newblock Second edition.
	
	\bibitem{FNP}
	E.~Feireisl, A.~Novotn{\' y}, and H.~Petzeltov{\' a}.
	\newblock On the existence of globally defined weak solutions to the
	{N}avier-{S}tokes equations of compressible isentropic fluids.
	\newblock {\em J. Math. Fluid Mech.}, {\bf 3}:358--392, 2001.
	
	\bibitem{FP20}
	E.~Feireisl and H.~Petzeltov\'{a}.
	\newblock On the long-time behaviour of solutions to the
	{N}avier-{S}tokes-{F}ourier system with a time-dependent driving force.
	\newblock {\em J. Dynam. Differential Equations}, {\bf 19}(3):685--707, 2007.
	
	\bibitem{FeiSwGw}
	E.~Feireisl and A.~\'Swierczewska-Gwiazda.
	\newblock The {R}ayleigh--{B}{\' e}nard problem for compressible fluid flows.
	\newblock {\em {\bf arxiv preprint No. 2110.10198}}, 2021.
	
	\bibitem{GAL}
	G.~P. Galdi.
	\newblock {\em An introduction to the mathematical theory of the {N}avier -
		{S}tokes equations, I.}
	\newblock Springer-Verlag, New York, 1994.
	
	\bibitem{GEHEHI}
	M.~Gei{\ss}ert, H.~Heck, and M.~Hieber.
	\newblock On the equation {${\rm div}\,u=g$} and {B}ogovski\u\i's operator in
	{S}obolev spaces of negative order.
	\newblock In {\em Partial differential equations and functional analysis},
	volume 168 of {\em Oper. Theory Adv. Appl.}, pages 113--121. Birkh\"auser,
	Basel, 2006.
	
	\bibitem{JinYang}
	Ch. Jin and T.~Yang.
	\newblock Time periodic solution to the compressible {N}avier-{S}tokes
	equations in a periodic domain.
	\newblock {\em Acta Math. Sci. Ser. B (Engl. Ed.)}, {\bf 36}(4):1015--1029,
	2016.
	
	\bibitem{KagOom}
	Y.~Kagei and R.~Oomachi.
	\newblock Stability of time periodic solution of the {N}avier-{S}tokes equation
	on the half-space under oscillatory moving boundary condition.
	\newblock {\em J. Differential Equations}, {\bf 261}(6):3366--3413, 2016.
	
	\bibitem{KagTsu}
	Y.~Kagei and K.~Tsuda.
	\newblock Existence and stability of time periodic solution to the compressible
	{N}avier-{S}tokes equation for time periodic external force with symmetry.
	\newblock {\em J. Differential Equations}, {\bf 258}(2):399--444, 2015.
	
	\bibitem{LI}
	P.-L. Lions.
	\newblock {\em Mathematical topics in fluid dynamics, Vol.1, Incompressible
		models}.
	\newblock Oxford Science Publication, Oxford, 1996.
	
	\bibitem{LI4}
	P.-L. Lions.
	\newblock {\em Mathematical topics in fluid dynamics, Vol.2, Compressible
		models}.
	\newblock Oxford Science Publication, Oxford, 1998.
	
	\bibitem{Tsuda}
	K.~Tsuda.
	\newblock Existence and stability of time periodic solution to the compressible
	{N}avier-{S}tokes-{K}orteweg system on {$\Bbb{R}^3$}.
	\newblock {\em J. Math. Fluid Mech.}, {\bf 18}(1):157--185, 2016.
	
	\bibitem{Valli1}
	A.~Valli.
	\newblock Navier-{S}tokes equations for compressible fluids: global estimates
	and periodic solutions.
	\newblock In {\em Nonlinear functional analysis and its applications, {P}art 2
		({B}erkeley, {C}alif., 1983)}, volume~45 of {\em Proc. Sympos. Pure Math.},
	pages 467--476. Amer. Math. Soc., Providence, RI, 1986.
	
	\bibitem{VAZA}
	A.~Valli and M.~Zajaczkowski.
	\newblock {N}avier-{S}tokes equations for compressible fluids: Global existence
	and qualitative properties of the solutions in the general case.
	\newblock {\em Commun. Math. Phys.}, {\bf 103}:259--296, 1986.
	
	\bibitem{WanWan}
	X.~Wang and W.~Wang.
	\newblock On global behavior of weak solutions to the {N}avier-{S}tokes
	equations of compressible fluid for {$\gamma=5/3$}.
	\newblock {\em Bound. Value Probl.}, pages 2015:176, 13, 2015.
	
	\bibitem{Ziemer}
	W.~Ziemer
	\newblock {\em Weakly differentiable functions}
	\newblock Springer-Verlag, 1989
	
\end{thebibliography}


\end{document}